# A generalization of the inverse tangent function

## 1) Background and motivation

In [1] we have proved the following formula for double integrals:

$$I(\alpha) = \int_0^\alpha \int_0^\alpha f(x,y)\,dx\,dy = \int_0^1 dx \int_0^\alpha \beta \{f(\beta,\beta x) + f(\beta x, \beta)\}\,d\beta \quad (F1)$$

valid for any function in two variables $f(x,y)$ continuous and bounded over a continuous domain containing all the squares $(0,\alpha) \times (0,\alpha)$, with $B > \alpha \geq 0$, $B$ being a positive real number. Similarly one can prove easily the more general formula

$$\int_0^\alpha \cdots \int_0^\alpha f(x_1,\cdots,x_n)\,dx_1\cdots dx_n = \int_0^1 \cdots \int_0^1 dx_1\cdots dx_n \int_0^\alpha \beta^{n-1}\{\Sigma\}\,d\beta \quad (F2)$$

where $\Sigma$ is the sum of the n functions:

$$\Phi_p = f(\beta x_1,\cdots,\beta x_{p-1},\beta,\beta x_{p+1},\cdots,\beta x_n)$$

Our aim is to use these two formulas to establish the value of some interesting sums of double or multiple integrals.

## 2) Generalized arctan functions

In this section we shall define a generalization of the usual *inverse tangent function* in one variable (usually denoted *arctan*): for any positive number $n \geq 2$, we define the $n^{th}$ order *arctan* as the $n - 1$ variables function:

$$\arctan(u_1, u_2, \cdots, u_{n-1}) = \int_0^{u_1} \cdots \int_0^{u_{n-1}} \frac{1}{1 + x_1^n + \cdots + x_{n-1}^n}\,dx_1 \cdots dx_{n-1}$$

This function is symmetric with respect to its variables.

To justify this definition, let us return briefly to the usual, single variable, *arctan* function, and more precisely to a new demonstration of one of its well known functional equations. The (F1) formula applied to the product of integrals on the left-hand side of the following relation yields (after taking the upper boundary $\alpha$ to infinity) :

$$2u\int_0^{+\infty} e^{-x^2}dx \int_0^{+\infty} e^{-(ux)^2}dx = \int_0^1 \frac{u}{u^2+x^2}dx + \int_0^1 \frac{u}{1+(ux)^2}dx = \arctan\frac{1}{u} + \arctan u$$

But we can compute each integral of the left-hand side separately, and find that its value is ½{Γ(½)}², that is $\pi/2$. Therefore, here we have proved the well known functional relation:

$$\arctan\frac{1}{u} + \arctan u = \frac{\pi}{2}$$

This proof is certainly more complicated than the classical one, as the latter involves only a simple derivation and the computation of the integration constant. But for us it is a discovery tool because it offers a mean to generalise this relation to a class of **several variables** functions where classical methods, such as differentition, do not work, and even fail to point to a clear direction for a generalization. In fact the main property involved here is the fundamental functional equation: $e^x . e^y = e^{x+y}$ of the exponential function, and indeed this can be generalised to the $n^{th}$ order *arctan* **basically because the exponential equation holds for any number of variables**: $e^a . e^b ...... e^z = e^{a+b+......+z}$.

For the sake of simplicity, and room on the paper, we shall treat completely only the case $n = 3$.

The (F2) formula allows us to write :

$$3uv \int_0^{+\infty} e^{-x^3}dx \int_0^{+\infty} e^{-(ux)^3}dx \int_0^{+\infty} e^{-(vx)^3}dx =$$

$$\int_0^1\int_0^1 \frac{uv}{1+(ux)^3+(vy)^3}dxdy + \int_0^1\int_0^1 \frac{uv}{x^3+u^3+(vy)^3}dxdy + \int_0^1\int_0^1 \frac{uv}{x^3+(uy)^3+v^3}dxdy$$

In each integral of the right-hand side we perform the changes of variables leading to the numerator $1 + t^3 + s^3$, and we obtain the sum:

$$\int_0^u \int_0^v \frac{1}{1+t^3+s^3} \, dt \, ds + \int_0^{\frac{1}{u}} \int_0^{\frac{v}{u}} \frac{1}{1+t^3+s^3} \, dt \, ds + \int_0^{\frac{1}{v}} \int_0^{\frac{u}{v}} \frac{1}{1+t^3+s^3} \, dt \, ds$$

that is:

$$\arctan(u,v) + \arctan(\frac{1}{u},\frac{v}{u}) + \arctan(\frac{1}{v},\frac{u}{v})$$

On the other hand, after the change of variables $x^3 = y$, we can see that:

$$3 \int_0^{+\infty} e^{-x^3} \, dx = \Gamma\left(\frac{1}{3}\right)$$

and compute the left-hand side of this relation as $3\{ ⅓ \, \Gamma(⅓) \}^3$. We are led to the formula:

$$\arctan(u,v) + \arctan(\frac{1}{u},\frac{v}{u}) + \arctan(\frac{1}{v},\frac{u}{v}) = 3 \left\{\frac{1}{3}\Gamma\left(\frac{1}{3}\right)\right\}^3$$

where the right-hand side is a constant playing the same role as $\pi/2$ in the classical formula.

For the general *arctan* function we have:

$$\arctan(u_1,u_2,\cdots,u_{n-1}) + \sum_{p=1}^{n-1} \arctan(\frac{u_1}{u_p},\frac{u_2}{u_p},\cdots,\frac{u_{p-1}}{u_p},\frac{1}{u_p},\frac{u_{p+1}}{u_p},\cdots,\frac{u_{n-1}}{u_p}) = n \left\{\frac{1}{n}\Gamma\left(\frac{1}{n}\right)\right\}^n$$

(demonstration left to the reader as an exercise).

If we set all the variables to 1 in this identity, it follows :

$$\int_0^1 \cdots \int_0^1 \frac{1}{1+x_1^n + \cdots + x_{n-1}^n} \, dx_1 \cdots dx_{n-1} = \left\{\frac{1}{n}\Gamma\left(\frac{1}{n}\right)\right\}^n$$

If we set all the variables equal to a given quantity $u$, and then make $u$ tends toward infinity, we cancel all the terms in the previous sum except the first one, and thus prove that:

$$\int_0^{+\infty} \cdots \int_0^{+\infty} \frac{1}{1+x_1^n + \cdots + x_{n-1}^n} \, dx_1 \cdots dx_{n-1} = n \left\{\frac{1}{n}\Gamma\left(\frac{1}{n}\right)\right\}^n$$

From the general functional relation satisfied by the $n^{th}$ order *arctan* function we can deduce solutions to certain functional equations.

For instance by setting $u_1 = 1/u$ and $u_p = 1$ for any $p > 1$, we find for the $n^{th}$ order *arctan* function:

$$\arctan(u, u, \cdots, u) + (n-1)\arctan(\frac{1}{u}, 1, \cdots, 1) = n\left\{\frac{1}{n}\Gamma\left(\frac{1}{n}\right)\right\}^n$$

For $u_1 = u$ and $u_p = 1$ for any $p > 1$:

$$\arctan(\frac{1}{u}, \frac{1}{u}, \cdots, \frac{1}{u}) + (n-1)\arctan(u, 1, \cdots, 1) = n\left\{\frac{1}{n}\Gamma\left(\frac{1}{n}\right)\right\}^n$$

It follows that if we define for any real positive *u*:

$$F(u) = \arctan(u, u, \cdots, u) + (n-1)\arctan(u, 1, \cdots, 1)$$

by addition of the two previous relations satisfied by the nth order *arctan* function we obtain:

$$F(u) + F(\frac{1}{u}) = 2n\left\{\frac{1}{n}\Gamma\left(\frac{1}{n}\right)\right\}^n$$

that is a class of solutions to the one-variable functional equation satisfied by the classical *arctan* function:

$$\Phi(x) + \Phi(\frac{1}{x}) = \text{constant}$$

If we consider the 4$^{th}$ order *arctan* function, we have for any positive real *u* and *v*:

$$\arctan(u, v, \frac{u}{v}) + \arctan(\frac{1}{u}, \frac{v}{u}, \frac{1}{v}) + \arctan(\frac{u}{v}, \frac{1}{v}, \frac{u}{v^2}) + \arctan(v, \frac{v^2}{u}, \frac{v}{u}) = 4\left\{\frac{1}{4}\Gamma\left(\frac{1}{4}\right)\right\}^4$$

If we define:

$$\arctan(u, v, \frac{u}{v}) + \arctan(\frac{u}{v}, \frac{1}{v}, \frac{u}{v^2}) = \Phi(u, v)$$

then, by taking into account the symmetry of the *arctan* function with regard to its variables, we find a two variables look alike of the classical second order *arctan* functional relation

$$\Phi(u, v) + \Phi(\frac{1}{u}, \frac{1}{v}) = 4\left\{\frac{1}{4}\Gamma\left(\frac{1}{4}\right)\right\}^4$$

*Juan PLA, 315 rue de Belleville, 75019 Paris (France)*